\documentclass[10pt,reqno]{amsart}
\usepackage{amssymb,amsmath,amscd}
\usepackage[english]{babel}

\newtheorem{lemma}{Lemma}
\newtheorem{theorem}{Theorem}
\newtheorem{definition}{Definition}

\newtheorem{proposition}{Proposition}
\def\id{\mathop{\rm id}}

\begin{document}

%\large

\title[Amenability of Closed Subgroups and Orlicz Spaces]
{Amenability of Closed Subgroups \\ and Orlicz Spaces}
\author{Yaroslav Kopylov}
\address{Sobolev Institute of Mathematics, Pr. Akad. Koptyuga 4,
630090, Novosibirsk, Russia and  Novosibirsk State University,
ul.~Pirogova 2, 630090, Novosibirsk, Russia }
\email{yakop@math.nsc.ru}
\thanks{The author was partially supported by the Russian Foundation for Basic Research 
(Grant~12-01-00873-a), the State Maintenance Program for the Leading Scientific Schools and
Junior Scientists of the Russian Federation (Grant~NSh-921.2012.1),
and the Integration Project 12-II-CO-01M-002 ``Geometric Analysis of Actual Problems 
in Function Theory and Differential Equations'' of the Siberian 
and Far Eastern Branches of the Russian Academy of Sciences.}

\date{}

\begin{abstract}
We prove that a closed subgroup $H$ of a second countable locally compact group~$G$ is amenable 
if and only if its left regular re\-pre\-sen\-ta\-tion on an Orlicz space $L^\Phi(G)$  for some 
$\Delta_2$-regular $N$-function $\Phi$ almost has invariant vectors. We also show 
that a noncompact second countable locally compact group~$G$ is amenable if and ony if 
the first cohomology space $H^1(G,L^\Phi(G))$ is non-Hausdorff for some $\Delta_2$-regular 
$N$-function $\Phi$. 
\end{abstract}

\maketitle

\section{Introduction}

Throughout, we assume all topological groups separated.

A locally compact topological group is called \textit{amenable}~\protect\cite{Eym}
if there exists a $G$-invariant mean on $L^\infty(G)$ or,
equivalently, $G$ possesses the \textit{fixed point property}:
for every action of~$G$ by continuous affine transformations on a nonempty convex compact 
subset~$Q$ of a locally convex space~$W$, there is a fixed point for $G$ in~$Q$.

Let $V$ be a Banach $G$-module, i.e., a real or complex Banach space
endowed with a continuous linear representation $\alpha:G\to \mathcal{B}(V)$.
We say that $V$ \textit{almost has invariant vectors} if, for every compact
subset $F\subset G$ and every $\varepsilon>0$, there exists a unit vector
$v\in V$ such that $\|\alpha(g)v-v\|\le\varepsilon$ for all $g\in F$.
Here $\mathcal{B}(V)$ stands for the space of all bounded linear
endomorphisms of a Banach space $V$.

Let $V$ be a~normed space of functions $f:G\to\mathbb{R}$ ($f:G\to\mathbb{C}$)
such that if $f\in V$ then, for every $g\in G$, the function
$$
\lambda_G(g)f(x)=f(g^{-1}x), \quad x\in G,
$$
lies in~$V$ and $\|\lambda_G(g)f\|_V = \|f\|_V$. Then $\lambda_G:G\to B(V)$
is called the {\it left regular representation} of~$G$ in~$V$.  

Examples of such function spaces $V$ are given by the space $L^p(G)$
of all real-valued functions on~$G$ integrable to the power
$p$ over $G$ with respect to a left-invariant Haar measure $\mu_G$. 
Instead of the $L^p$ spaces, one can consider more general 
Orlicz spaces $L^\Phi(G)$ of real-valued functions on~$G$ with the finite ``gauge'' norm
$$
\|f\|_{(\Phi)}= \inf \biggl\{ k>0 \, : 
\, \int\limits_G \Phi\left(\frac{f(x)}{k}\right) d\mu_G(x) \le 1 \biggr\}.
$$
for an~$N$-function $\Phi$ ($N$-functions are defined in Section~\ref{N-func}). 
Orlicz spaces on locally compact groups were
considered in~\cite{Bu75,KamMus89} and more recently in~\cite{Ak2012,Ak2013,Ra01,Ra04a}.
 
In~\cite{St65} (see also \cite[Theorem~8.3.2]{ReiSt}), Stegeman proved that, 
for a~locally compact group, the following conditions $(P_p)$ (called Reiter's conditions)
are equivalent for all~$p\ge 1$:

($P_p$) for every compact set $F\subset G$ and every $\varepsilon>0$, there
exists a function $f\in L^p(G)$ with $f\ge 0$ and $\|f\|_{L^p(G)}=1$
such that $\|\lambda_G(z)f-f\|_{L^p(G)}<\varepsilon$ for all $z\in F$.

In~\cite{Eym} Eymard extended this equivalence to quotients $G/H$ of locally
compact groups by closed subgroups and proved that conditions~($P_p$) are
equivalent to the amenability of~$G/H$.

In~\cite[Proposition~2, pp.~387--389]{Ra04a}, Rao proved that a locally compact 
group $G$ is amenable if and only if, given a $\Delta_2$-regular 
$N$-function $\Phi$, \,$G$ satisfies the property 

($P_\Phi$) for every compact set $F\subset G$ and every $\varepsilon>0$, there
exists a function $f\in L^\Phi(G)$ with $f\ge 0$ and $\|f\|_{L^{(\Phi)}(G)}=1$
such that $\|\lambda_G(z)f-f\|_{L^{(\Phi)}(G)}<\varepsilon$ for all $z\in F$.

Here $\|\cdot\|_{L^{(\Phi)}(G)}$ stands for the gauge norm in the space $L^\Phi(G)$.
 
In~2005, Bourdon, Martin, and Valette established the following \protect\cite[Lemma~2]{BMV}:

\smallskip

{\bf Theorem A.}
{\it Suppose that $p\ge 1$. Let $X$ be a countable set on
which a countable group $H$ acts freely. The following are equivalent:

{\rm(i)} The natural ``permutation'' representation $\lambda_X$ of $H$
on $L^p(X)$ almost has invariant vectors;

{\rm(ii)} $H$ is amenable. }
\smallskip

In~\cite{Kop05}, in an attempt to generalize this assertion, we proved:

\smallskip 

{\bf Theorem B.}
Assume that $p\ge 1$. Let $G$ be a second countable locally compact group
and let $H$  be a closed subgroup in $G$. The following are equivalent:

{\rm(i)} The left regular representation of $H$ on $L^p(G)$ almost has invariant vectors;

{\rm(ii)} $H$ is amenable.

\smallskip

The paper is organized as follows: In Section~\ref{N-func}, we recall some basic notions 
concerning $N$-functions and Orlicz 
spaces. Section~\ref{integr} contains some necessary in\-for\-ma\-tion on integration on locally 
compact groups and homogeneous spaces. In Section~\ref{amen-subgr}, we prove a generalization 
of Theorem~B, where $L^p(G)$ is replaced by the Orlicz space $L^\Phi(G)$ for any 
$\Delta_2$-regular $N$-function $\Phi$ (Theorem~\ref{amenab-subgr}). In Section~\ref{cohom-amen},  
using the equivalence of amenability and the fulfillment of the above condition~($P_\Phi$) and 
a general result by Guichardet, we deduce that the nonreduced and reduced first cohomology 
of a noncompact second countable locally compact group coincide if and only if it is not amenable.

\section{$N$-Functions and Orlicz Function Spaces}\label{N-func}

\begin{definition}
A function $\Phi:\mathbb{R}\to \mathbb{R}$ is called an {\it $N$-function} if

{\rm(i)} $\Phi$ is even and convex; 

{\rm(ii)} $\Phi(x)=0 \Longleftrightarrow x=0$; 

{\rm(iii)} $\lim\limits_{x\to 0} \frac{\Phi(x)}{x}=0; 
\quad \lim\limits_{x\to \infty} \frac{\Phi(x)}{x}=\infty$.
\end{definition}

An $N$-function $\Phi$ has left and right derivatives (which can differ only on an at most
countable set, see, for instance, \cite[Theorem~1, p.~7]{RaRen91}).
The left derivative $\varphi$ of $\Phi$ (we write $\varphi=\Phi'$ below) is left continuous, 
nondecreasing on~$(0,\infty)$, and such that $0<\varphi(t)< \infty$ for $t>0$, $\varphi(0)=0$,
$\lim\limits_{t\to\infty}\varphi(t)=\infty$. The function 
$$
\psi(s) = \inf \{t>0 \,:\, \varphi(t)>s \}, \quad s>0,
$$
is called the {\it left inverse} of~$\varphi$. 

The functions $\Phi, \Psi$ given by
$$
\Phi(x) = \int\limits_0^{|x|} \varphi(t) dt, \quad \Psi(x) = \int\limits_0^{|x|} 
\psi(t) dt
$$
are called {\it complementary $N$-functions}.

The $N$-function $\Psi$ complementary to an $N$-function $\Phi$ can also be expressed as
$$
\Psi(y)= \sup \{ x|y| - \Phi(x) \,:\, x\ge 0\}, \quad y\in\mathbb{R}.
$$

$N$-functions are classified in accordance with their growth rates as follows:

\begin{definition}
An $N$-function $\Phi$ is said to satisfy the $\Delta_2$-condition for large~$x$
(for small~$x$, for all~$x$), which is written as $\Phi\in\Delta_2(\infty)$ 
($\Phi\in\Delta_2(0)$, or $\Phi\in\Delta_2$), if there exist constants $x_0>0$,
$K>2$ such that $\Phi(2x)\le K\varphi(x)$ for $x\ge x_0$ (for $0\le x \le x_0$,
or for all~$x\ge 0$); and it satisfies the $\nabla_2$-condition for large~$x$ (for small~$x$, 
or for all~$x$), denoted symbolically as $\Phi\in\nabla_2(\infty)$ 
($\Phi\in\nabla_2(0)$, or $\Phi\in\nabla_2$) if there are constants $x_0>0$ and
$c>1$ such that $\Phi(x)\le \frac{1}{2c}\Phi(cx)$ for $x\ge x_0$ (for $0\le x\le x_0$,
or for all~$x\ge0$).
\end{definition}

\medskip

Henceforth, let $\Phi$ be an $N$-function and let $(\Omega,\Sigma, \mu)$ be a measure
space. 

\begin{definition}
The set $\tilde L^\Phi=\tilde L^\Phi(\Omega)= L^\Phi(\Omega,\Sigma,\mu)$ is defined 
to be the set of measurable functions $f:\Omega \to \mathbb{R}$ such that
$$
\rho_{\Phi}(f):= \int\limits_\Omega \Phi(f) d\mu <\infty.
$$
\end{definition}

\begin{proposition}~\cite{RaRen02}  
The set $\tilde L^\Phi$ is a vector space in the following cases:

{\rm(i)} $\mu(\Omega)< \infty, \quad \Phi\in \Delta_2(\infty)$;

{\rm(ii)} $\mu(\Omega)= \infty, \quad \Phi\in \Delta_2$. 

{\rm(iii)} $\Omega$ is countable, $\mu$ is the~counting measure on~$\Omega$, \quad $\Phi\in \Delta_2(0)$. 

\end{proposition}

\begin{definition}
The linear space 
$$
L^\Phi= L^\Phi(\Omega) = L^\Phi(\Omega,\Sigma,\mu) = \{ f:\Omega\to\mathbb{R} ~\text{measurable}~ : 
\rho_\Phi(af)<\infty ~\text{for {\it some} $a>0$}\} 
$$
is called an {\it Orlicz space} on~$(\Omega,\Sigma,\mu)$. 
\end{definition}

For an Orlicz space $L^\Phi=L^\Phi(\Omega,\Sigma,\mu)$, the $N$-function $\Phi$ is called 
{\it $\Delta_2$-regular} if $\Phi\in\Delta_2(\infty)$ when $\mu(\Omega)<\infty$ or                        
$\Phi\in\Delta_2$ when $\mu(\Omega)=\infty$ or $\Phi\in\Delta_2(0)$ for $\mu$ 
a~counting measure.

Let $\Psi$ be the complementary $N$-function to~$\Phi$.  

Below we as usual identify two functions equal on a~set of measure zero.

If $f\in L^\Phi$ then the functional $\|\cdot\|_\Phi$ (called 
{\it the Orlicz norm}) defined by 
$$
\|f\|_\Phi= \|f\|_{L^\Phi(\Omega)} = \sup \biggl\{ \biggl| \int\limits_\Omega fg \, d\mu \biggr| : 
\rho_\Psi(g) \le 1 \biggr\}
$$
is a seminorm. It becomes a norm if $\mu$ satisfies the {\it finite subset property}
(see~\cite[p.~59]{RaRen91}): if $A\in\Sigma$ and $\mu(A)>0$ then there exists $B\in\Sigma$,
$B\subset A$, such that $0<\mu(B)<\infty$. 

The {\it gauge} (or {\it Luxemburg}) {\it norm} of a function $f\in L^\Phi$
is defined by the formula
$$
\|f\|_{(\Phi)}= \|f\|_{L^{(\Phi)}(\Omega)} =
\inf \biggl\{ k>0 \, : \, \rho_\Phi\biggl(\frac{f}{k}\biggr)\le 1 \biggr\}.
$$
This is a norm without any constraint on the measure~$\mu$ (see \cite[p.~54, Theorem~3]{RaRen91}). 

Suppose that the measure~$\mu$ satisfies the finite subset property.
As is proved in~\cite[Chapter~10]{Ra04b}, a left-invariant Haar measure on a~locally compact group 
has this property.

It is well known that the Orlicz and gauge norms are equivalent, namely (see, for example,
\cite[pp.~61--62]{RaRen91}):
$$
\|f\|_{(\Phi)} \le \|f\|_\Phi  \le  2 \|f\|_{(\Phi)}.
$$

We will need the following version of H\"older's inequality for Orlicz 
spaces~\cite[p.~62]{RaRen91}:

\smallskip

{\bf H\"older's Inequality.}  
{\it If $\Phi$ and $\Psi$ are two complementary $N$-functions then $fg\in L^1$ and
\begin{equation}\label{holder}
\|fg\|_1 \le \|f\|_{(\Phi)} \|g\|_\Psi \quad (~\|fg\|_1 \le \|f\|_\Phi \|g\|_{(\Psi)}~). 
\end{equation}                      
}

\section{Integration on Locally Compact Groups and Borel Sections}\label{integr}

Recall some basic facts and definitions from the theory of integration 
on locally compact groups.

Let $G$ be a locally compact group and let $H$ be a closed subgroup
in~$G$. Denote by $\mu_G$ and $\mu_H$ left-invariant Haar measures
on~$G$ and~$H$ respectively and denote by~$\pi$ the projection
$G\to G/H$.

Denote by~$\Delta_K$ the modulus of a locally compact group~$K$.

Given a function $f$ and a class $u\in G/H$, take an arbitrary representative $x$ in $u$
and consider the function $\alpha:y\to f(xy)$ on~$H$. If $\alpha$
is integrable over~$H$, the left invariance of~$\mu_H$ implies that
$\int\limits_H f(xy)\,d\mu_H(y)$ is independent of the choice of~$x$ with
$\pi(x)=u$.

It is well known that the homogeneous space $G/H$ admits a {\it quasi-$G$-invariant}
measure $\mu_{G/H}$ on~$H$ which is unique up to equivalence.
Here the ``quasi-$G$-invariance'' means that all left translates of~$\mu_{G/H}$
by the elements of~$G$ are equivalent to~$\mu_{G/H}$. The measure $\mu_{G/H}$
can be described as follows (see~\cite[Chapter~VII, 2.5]{B2} or~\cite{Eym}).

\smallskip
(a) There exists a positive continuous function $\rho$ on~$G$ such that
$\rho(xy)=\frac{\Delta_H(y)}{\Delta_G(y)}{\rho(x)}$
for all $x\in G$ and $y\in H$.
\smallskip

Put $\mu_{G/H}=(\rho\mu_G)/\mu_H$ (see~\cite[Definition~1 in Chapter~VII, 2.2]{B2}).

\smallskip
(b) If $f\in L^1(G,\rho\mu_G)$ then the set of $\overline{x}=\pi(x)\in G/H$
for which $y\mapsto f(xy)$ is not $\mu_H$-integrable is
$\mu_{G/H}$-negligible, the function
$\overline{x}=\pi(x)\mapsto \int\limits_H f(xy)d\mu_H(y)$
is $\mu_{G/H}$-integrable, and
$$
\int\limits_G f(x)\rho(x) \,d\mu_G(x)
= \int\limits_{G/H} \,d\mu_{G/H}(\overline{x})\int\limits_H f(xy) \,d\mu_H(y).
$$
\smallskip

(c) There exists a nonnegative continuous function $h$ on~$G$ with
$\int\limits_H h(xy)dy=~1$ for all $x\in G$ such that a function
$w$ on~$G/H$ is $\mu_{G/H}$-measurable ($\mu_{G/H}$-integrable)
if and only if $h(w\circ\pi)$ is $\rho\mu_G$-measurable
($\rho\mu_G$-integrable). If
$w\in L^1(G/H,\mu_{G/H})$ then
$$
\int\limits_{G/H} w(\bar{x}) \,d\mu_{G/H}(\bar{x})
= \int\limits_G h(x)w(\pi(x))\rho(x) \,d\mu_G(x).
$$

Note that a second countable locally compact space is Polish (\textit{polonais})
(see~\cite{B2}). As follows from Dixmier's lemma (see~\cite{Di62}), if $G$
is a Polish group and $H$ is a closed subgroup
in~$G$ then there exists a Borel section $\sigma:G/H\to G$
(in particular, $\pi\circ\sigma=\id_{G/H}$). We will need the following
technical assertion (see~\cite{Kop05} for a proof):

\begin{lemma}\label{l1}
Suppose that $G$ is a second countable locally compact group,
$H$ is a closed subgroup in~$G$, $\sigma:G/H\to G$ is a Borel
section, and $f\in L^1(G,\rho\mu_G)$. Then, in the above notations,
$$
\int\limits_G f(x)\rho(x) \,d\mu_G(x) = \int\limits_H d\mu_H(y)
\int\limits_{G/H} f(\sigma(\overline x) y) \,d\mu_{G/H}(\overline x).
$$
\end{lemma}

\section{Amenability of Closed Subgroups}\label{amen-subgr}

The main result of this section is as follows:

\begin{theorem}\label{amenab-subgr}
Assume that $\Phi$ is a $\Delta_2$-regular $N$-function.
Let $G$ be a second countable locally compact group
and let $H$ be a closed subgroup in $G$. The following are equivalent:

{\rm(i)} The left regular representation of $H$ on $L^\Phi(G)$ almost has invariant vectors;

{\rm(ii)} $H$ is amenable.
\end{theorem}

\textit{Proof.}
Put $\Phi'=\varphi$ and let $\Psi$ be the complementary $N$-function to~$\Phi$.

Observe first that (ii) implies (i) by the equivalence of amenability and the fulfillment of 
the Rao--Reiter condition~($P_\Phi$), established by Rao 
in~\cite[Proposition~2, pp.~387--389]{Ra04a}:  

($P_\Phi$) For every compact set $F$ and every $\varepsilon>0$, there
exists a function $f\in L^\Phi(H)$ with $f\ge 0$ and $\|f\|_{L^{(\Phi)}(H)}=1$
such that $\|\lambda_H(z)f-f\|_{L^{(\Phi)}(H)}<\varepsilon$ for all
$z\in F$.

Now, prove (i)$\Rightarrow$(ii). Suppose that $L^\Phi(G)$ almost
has invariant vectors for $H$ and deduce from this that $H$ meets
Reiter's condition~($P_1$).

By the $\Delta_2$-regularity of~$\Phi$, we conclude from 
\cite[Proposition 8, p.~79]{RaRen91} that 
$$
S:= \sup \bigl\{ \rho_\Psi(\varphi\circ |v|) \, : \, v\in L^\Phi(G), 
\, \|v\|_{L^{(\phi)}(G)}\le 1 \bigr\} <\infty.
$$ 

Take $\varepsilon>0$ and a compact set $F\subset H$; choose
$f\in L^\Phi(G)$, $\|f\|_{L^{(\Phi)}(G)}=1$, such that
\begin{equation}\label{ff}
\|\lambda_G(z)f-f\|_{L^{(\Phi)}(G)}\le\frac{\varepsilon}{2(S+1)}
\end{equation}
for all $z\in F$. Assume without loss of generality that $f\ge 0$ (taking
$|f|$ instead of~$f$ if necessary). 

Note that 
\begin{equation}\label{ineqq}
|\Phi(a)-\Phi(b)| \le |a-b|\,(\varphi(a)+\varphi(b)), \quad a,b\ge 0, 
\end{equation}
because $\varphi$ is monotone and nonnegative (cf. \cite[p. 388]{Ra04a}). 

Put $u=\Phi\circ f$.  Since $\Phi$ is $\Delta_2$-regular, we have 
(\cite[p.~78]{KraRu}):  
$$
\int\limits_G u(x) \,d\mu_G(x) = \int\limits_G \Phi(f(x)) \,d\mu_G(x) = \|f\|_{L^{(\Phi)}(G)}= 1.
$$               
For $z\in F$, using~(\ref{ineqq}),
H\"older's inequality~(\ref{holder}), (\ref{ff}), and the inequality 
$$
\|v\|_{L^{\Psi}(G)} \le \rho_\Psi(v)+1,
$$
we infer
\begin{multline*}
\|\lambda_G(z)u-u\|_{L^1(G)}
= \int\limits_G |\Phi(f(z^{-1}x))-\Phi(f(x))| \,d\mu_G(x)
\\
\le  \int\limits_G |f(z^{-1} x)-f(x)|\,\, |\varphi(f(z^{-1}x))+\varphi(f(x))| \,\, d\mu_G(x)
\\
\le \|\lambda_G(z)f-f\|_{L^{(\Phi)}(G)} \|\varphi\circ \lambda_G(z)f + \varphi\circ f\|_{L^\Psi(G)}
\\
\le \|\lambda_G(z)f-f\|_{L^{(\Phi)}(G)}
\bigl( \|\varphi\circ \lambda_G(z)f\|_{L^\Psi(G)} +\|\varphi\circ f\|_{L^\Psi(G)} \bigr)
\\
= 2\|\varphi\circ f\|_{L^\Psi(G)} \|\lambda_G(z)f-f\|_{L^{(\Phi)}(G)}
\le 2(\rho_\psi(\varphi\circ f) + 1) \, \|\lambda_G(z)f-f\|_{L^{(\Phi)}(G)} 
\\
\le 2(S+1) \|\lambda_G(z)f-f\|_{L^{(\Phi)}(G)} < \varepsilon.
\end{multline*}

Now, let $\sigma:G/H\to G$ be a Borel section. Consider the function 
$$
U(y)= \int\limits_{G/H} \frac{u(y\sigma(\overline x))}{\rho(\sigma(\overline x))} 
\,d\mu_{G/H}(\overline x), \quad y\in H,
$$
where $\rho$ is the function described in Section~\ref{integr}. By Lemma~\ref{l1}, since
$u$ is nonnegative, we have 
$$
\|U\|_{L^1(H)}=\int\limits_G u(x) d\mu_G(x) = 1.
$$
Involving Lemma~\ref{l1} again, we obtain the following estimates:
\begin{multline*}
\|\lambda_H(z)U-U||_{L^1(H)}
= \int\limits_H  \biggl| \int\limits_{G/H}
\frac{u(z^{-1}y\sigma(\overline x))-u(y\sigma(\overline x))}{\rho(\sigma(\overline x))}
\,d\mu_{G/H}(\overline x) \biggr| \,d\mu_H(y)
\\
\le \int_H \limits d\mu_H(y) \int\limits_{G/H}  \biggl|
\frac{u(z^{-1}y\sigma(\overline x))-u(y\sigma(\overline x))}{\rho(\sigma(\overline x))} \biggr|
\,d\mu_{G/H}(\overline x)
\\
= \int\limits_G \frac{|u(z^{-1}x)-u(x)|}{\rho(x)}\rho(x) \,d\mu_G(x)
= \int\limits_G |u(z^{-1}x)-u(x)| \,d\mu_G(x)
= \|\lambda_H(z)u-u\|_{L^1(G)}.
\end{multline*}
So, if $z\in F$ then $\|\lambda_H(z)U-U\|_{L^1(H)}<\varepsilon$.
Thus, $H$ has property~($P_1$) and hence is amenable.
Theorem~\ref{amenab-subgr} is proved.
\qed

{\bf Remark.}
Theorem~\ref{amenab-subgr} is informative only if $H$ is noncompact since 
(i) and~(ii) are both fulfilled when $H$ is compact.

\section{First Cohomology and Amenability}\label{cohom-amen}

Let $G$ be a topological group and let~$V$ be a~topological $G$-module, i.e., a real or complex
topological vector space endowed with a~linear representation $\pi:G\times V\to V$, 
$(g,v)\mapsto \pi(g)v$. The space~$V$ is called a {\it Banach $G$-module} if $V$ is a~Banach 
space and $\pi$ is a~representation of~$G$ by isometries of~$V$. Introduce the notation:
$$
Z^1(G,V):=\{b:G\to V \text{~continuous~} \mid b(gh)=b(g) + \pi(g)b(h) \} \quad 
\text{({\it 1-cocycles})}; 
$$
$$
B^1(G,V) = \{ b\in Z^1(G,V) \mid (\exists v\in V)\, (\forall g\in G) \,\, b(g)= \pi(g) v \}
\quad \text{({\it 1-coboundaries})};
$$
$$
H^1(G,V) = Z^1(G,V)/B^1(G,V) \quad \text{({\it 1-cohomology with coefficients in~$V$})}.
$$

Endow $Z^1(G,V)$ with the topology of uniform convergence on compact subsets of~$G$ and 
denote by $\overline{B}^1(G,V)$ the closure of $B^1(G,V)$ in this topology. The quotient 
$\overline{H}^1(G,V) = Z^1(G,V)/\overline{B}^1(G,V)$ is called the 
{\it reduced 1-cohomology} of~$G$ with coefficients in the $G$-module~$V$. 

The following assertion was established by Guichardet (see~\cite[Th\'eor\`eme~1]{Gui72}):

\begin{lemma}\label{guich-th}
Let $G$ be a~locally compact second countable group and let $V$ be a~Banach module 
such that
$$
V^G := \{ v\in V \mid \pi(g)v=v \text{~for all~} g\in G \} =0.
$$
Then the following are equivalent:

{\rm(i)} $H^1(G,V) = \overline{H}^1(G,V)$;

{\rm(ii)} $V$ does not almost have invariant vectors, that is, there exists a~compact 
subset~$F\subset G$ and $\varepsilon>0$ such that 
$\sup_{g\in F} \|\pi(g)v-v\| \ge \varepsilon\|v\|$ for all $v\in V$. 
\end{lemma}

As is well known, if a locally compact group~$G$ is noncompact then the Haar measure 
of the whole group is infinite~\cite[Theorem~15.9]{HRI}. Hence, constant functions
are not integrable over~$G$. Therefore, $L^\Phi(G)^G=0$. Thus, combining 
Lemma~\ref{guich-th} with the Rao--Reiter condition~($P_\Phi$),  
we obtain the following generalization of~Corollary~2.4 in \cite[p.~86]{MV07} to 
coefficients in an~Orlicz space:

\begin{proposition}\label{or-cohom}
Suppose that $\Phi$ is a $\Delta_2$-regular $N$-function.
If $G$ is a noncompact second countable locally compact group then the following 
are equivalent:
           
{\rm(i)} $H^1(G,L^\Phi(G)) = \overline{H}^1(G,L^\Phi(G))$;

{\rm(ii)} $G$ is not amenable.
\end{proposition}

\smallskip

The author is indebted to the referee for valuable remarks.

\end{document}